\documentstyle{amsppt}
\magnification =1200
\title Level Sets and the Uniqueness of Measures
\endtitle

\author Dale E. Alspach
\endauthor
\thanks
Research supported in part 
by NSF grant DMS-890327
\endthanks
\address
Department of Mathematics
Oklahoma State University
Stillwater, OK  74078-0613
\endaddress
\abstract
A result of Nymann is extended to 
show that a positive $\sigma$-finite measure with
range an interval is determined by its level sets.  
An example is given of two finite positive
measures with range the same finite union of intervals but 
with the property that one is
determined by its level sets and the other is not.
\endabstract
\endtopmatter
\document

\heading {\bf 0.  Introduction}\endheading
The purpose of this note is to extend some results of
Leth [L], Malitz [M], and Nymann [N] and answer some questions
raised in these papers.  Those papers deal with a uniqueness
property exhibited by positive measures whose range is an interval
and used techniques from real analysis.  The main idea in this paper
to change from a purely real analysis approach to the problems to a
functional analytic one.  This approach reveals and clarifies the
issues involved.
 
The basic problem is the following.  Suppose that $\mu$ is a $\sigma$-finite
signed measure on some measurable space $(\Omega ,\Cal B)$.  We will say
that a signed measure $\nu$ on $(\Omega ,\Cal B)$ satisfies $(\Cal L)$ if for every $
A$ and
$B\in \Cal B$ such that $\mu (A)=\mu (B)\neq\pm\infty ,$ $\nu (A)=
\nu (B)\neq\pm\infty .$

\remark{Question 1}  Under what assumptions on $\mu$ is it true that if $\nu$
satisfies $(\Cal L)$, $\nu =\alpha\mu$ 
for some $\alpha\in \Bbb R$?
\endremark 
\medskip

We will also consider the stronger condition $(\Cal O)$ which requires that
condition $(\Cal L)$ hold and for every $A$ and $B\in \Cal B$ such that
$\mu (A)\leq \mu (B)\text{, }\nu (A)\leq \nu (B).$

\remark{Question 2}  Under what assumptions on $\mu$ is it true that if $\nu$
satisfies $(\Cal O)$, 
$\nu =\alpha\mu$ for some $\alpha \in \Bbb R$?
\endremark
\medskip

As a convenience in stating results we will say that $\mu$ is {\it uniquely
determined by $(\Cal L)$\/}, respectively, $(\Cal O)$, if $\mu$ satisfies some property P
and for any $\nu$ satisfying $(\Cal L)$, respectively, $(\Cal O)$, $\nu$ is 
a constant times
$\mu .$

The condition $(\Cal L)$ was used in [N] where it was shown that for finite
positive measures with range equal to an interval that $\mu$ was
uniquely determined by $(\Cal L)$.  The symmetric form of $(\Cal O)$ was
considered by Leth for the case of positive purely atomic measures
such that for each $\varepsilon >0$, there are only finitely many atoms of mass
greater than $\varepsilon .$ He called measures satisfying the symmetric form of
$(\Cal O)$ sympathetic.  This condition arose from some problems in
qualitative measure theory.  (See the references in [L] for more on
this.)  Malitz [M] showed that for this same class of atomic
measures the symmetric form of $(\Cal L)$ implies the symmetric form of
$(\Cal O)$.

It is easy to give examples where even under the symmetric form of
$(\Cal O)$ $\mu$ is not uniquely determined.  As pointed about by Leth [L], if
$0<r<\frac 12$ and $\mu$ is purely atomic with atoms $\{A_n:n=0,1,
...\}$ such that
$\mu (A_n)=r^n,$ there are no sets $A\neq B$ such that $\mu (A)=\mu 
(B)<\infty .$  In
order to generalize this observation, Leth makes the following
definition.  An atom $A$ is a {\it bully\/} if
$$\sup \{\sum\mu (B_j):(B_j)\subset \Cal B,B_j\cap B_k=\emptyset ,
j\neq k\text{ and }\mu (B_j)<\mu (A)\text{ for all }j\}<\mu (A).$$
In the example just given every atom is a bully.  Clearly if $\mu$ is
purely atomic and every atom is a bully, the condition $(\Cal L)$ in vacuous
and any purely atomic measure $\nu$ for which every atom is a bully
will satisfy $(\Cal O)$.  On the other hand, Leth gives examples to show
that bullies can be present but under $(\Cal L)$ $\mu$ is uniquely determined.
Leth (as improved by Malitz) also showed that if the space consists
of a sequence of atoms with measures converging to 0 and there are
no bullies then $\mu$ is uniquely determined by the symmetric
form of $(\Cal L)$.  It is not hard to see
that the no bullies condition is equivalent for positive measures to
the condition that the range of the measure is an interval and
therefore Nymann's result can be viewed as an extension of Leth's in
the case of finite measures.  Thus the question remains as to what
conditions on $\mu$ imply that $(\Cal L)$ or $(\Cal O)$ uniquely determines $
\mu$.

In this paper we will complete the work of [N] by showing that a
positive $\sigma$-finite measure with range an interval is uniquely
determined by $(\Cal L)$.  We will also prove some results on signed
measures.  Leth asked if it is possible to construct a measure so
that if the atoms are listed in decreasing order every second atom is
a bully, but the measure is uniquely determined.  We will show that
this is impossible.  We also will give several examples which show
that the range of measure is not a good indicator of whether the
measure is uniquely determined.  In particular we present an example
of two positive measures $\mu$ and $\nu$ such that the range of $\mu$
is a finite union
of intervals, the range of $\mu$ equals the range of $\nu ,$ $\mu$ is uniquely
determined by $(\Cal L)$ but $\nu$ is not uniquely determined by $(\Cal O)$.

We will use standard notation and terminology from real analysis and
elementary functional analysis as may be found in [F] or [R].  By
a $\sigma$-finite signed measure $\mu$ we mean an extended real valued
measure on a measurable space $(\Omega,\Cal B)$ taking on only one
of the two infinite values $\pm\infty$ such that
$|\mu |$ is $\sigma$-finite.  If $\mu$ is a $\sigma$-finite signed 
measure we will denote
the atoms of $\mu$ by $\{A_i\}$ and the complement of the union of the atoms
by $C$.  If for every $\varepsilon >0$ there only finitely many atoms $
A_i$ with
$|\mu |(A_i)>\varepsilon ,$ we will assume that the atoms are ordered so that
$|\mu |(A_i)\geq |\mu |(A_{i+1})$ for all i.  
$L_p$ will always refer to
$L_p(\Omega ,\Cal B,|\mu |),1\leq p\leq \infty .$
If $R$ is a subset of a vector space $X$ then sp $R$ will
denote  the span of $R$.

\newpage
\heading {\bf 1.  A Functional Analytic Approach.}\endheading

Before we reformulate the problem in functional analytic terms,
we will present two simple lemmas which will allow us to assume
that the measure $\mu$ is positive in most of our considerations.
We would like to thank the referee who suggested these to us and
made several other suggestions to improve the exposition.

\proclaim{Lemma 1.1} Suppose that $\mu$ is a $\sigma$-finite
signed measure on $(\Omega,\Cal B)$ and $\Omega^+$ and
$\Omega^-$ are the positive and negative sets of the Hahn
decomposition of $\mu$, i.e., $\mu_{|\Omega^+}\geq 0$,
$\mu_{|\Omega^-}\leq 0$ and $\Omega^+\cup\Omega^-=\Omega$. Suppose that $\nu$ is another
$\sigma$-finite
signed measure on $(\Omega,\Cal B)$ and let
$$ \nu'(A)=\nu(A\cap \Omega^+)-\nu(A\cap \Omega^-)$$
for all $A\in \Cal B.$ Then $\nu'$ is a $\sigma$-finite
signed measure on $(\Omega,\Cal B)$ such that $\nu$ satisfies $(\Cal L)$,
repsectively, $(\Cal O)$,
with respect to $\mu$ if and only if $\nu'$ satisfies $(\Cal L)$,
respectively, $(\Cal O)$, with respect to $|\mu|$.
\endproclaim
\demo{Proof} Suppose that $(\Cal O)$ is satisfied for $\mu$ and
$\nu$ and that $A, B \in \Cal B$ such that 
$$|\mu|(A)=\mu(A\cap \Omega^+)-\mu(A\cap \Omega^-)\leq\mu(B\cap
\Omega^+)-\mu(B\cap \Omega^-)=|\mu|(B).$$
Then
$$\mu(A\cap \Omega^+)+\mu(B\cap \Omega^-)\leq\mu(B\cap
\Omega^+)+\mu(A\cap \Omega^-).$$
By $(\Cal O)$,
$$\nu(A\cap \Omega^+)+\nu(B\cap \Omega^-)\leq\nu(B\cap 
\Omega^+)+\nu(A\cap \Omega^-)$$
and hence
$$\nu'(A)=\nu(A\cap \Omega^+)-\nu(A\cap \Omega^-)\leq\nu(B\cap
\Omega^+)-\nu(B\cap \Omega^-)=\nu'(B).$$
The proofs of the converse and the case of $(\Cal L)$ are
similar.
\qed\enddemo

Thus $\mu$ is uniquely determined by $(\Cal L)$, respectively,
$(\Cal O)$, if and only if $|\mu|$  is uniquely determined by
$(\Cal L)$, respectively,
$(\Cal O)$, 

Most of the conditons we will impose on $\mu$ are on the range of
$\mu$ and the next lemma shows that the conditions are easily
transferred to the range of $|\mu|$.

\proclaim{Lemma 1.2} If $\mu$ is a $\sigma$-finite
signed measure on $(\Omega,\Cal B)$ and $\Omega^+$ and
$\Omega^-$ are the positive and negative sets of the Hahn
decomposition of $\mu$, then 
$$\alignat2
\text{range }|\mu|&= \text{range }
\mu-\mu(\Omega^-),&&\qquad\text{if } \mu(\Omega^-)>-\infty,
\\
 \intertext{and}
\text{range }|\mu|&= -\text{range }
\mu+\mu(\Omega^+),&&\qquad\text{if }\mu(\Omega^+)<+\infty.
\endalignat$$
\endproclaim
\demo{Proof} Assume that $\mu(\Omega^-)>-\infty.$ If $A \in \Cal
B,$
$$ \align
\mu(A)=\mu(A\cap \Omega^+)+\mu(A\cap \Omega^-)
&=\mu(A\cap \Omega^+)+\mu(\Omega^-)-\mu( \Omega^-\setminus A)\\
&=|\mu|((A\cap \Omega^+)\cup (\Omega^-\setminus A))+\mu(\Omega^-)
\endalign$$
and
$$ \align
|\mu|(A)=\mu(A\cap \Omega^+)-\mu(A\cap \Omega^-)
&=\mu(A\cap \Omega^+)-\mu(\Omega^-)+\mu( \Omega^-\setminus A)\\
&=\mu((A\cap \Omega^+)\cup (\Omega^-\setminus A))-\mu(\Omega^-) .
\endalign$$
Thus $\text{range }|\mu|= \text{range }
\mu-\mu(\Omega^-)$. The other case is similar.
\qed\enddemo

We will now formulate the problem in functional analytic terms.  
Let $\mu$ be a  positive $\sigma$-finite 
measure and
 consider the 
following linear subspace X of
$L_0$, the space of equivalence classes (a.e. $\mu $) of measurable functions
on $(\Omega ,\Cal B).$
$$X=\text{sp}\{1_A-1_B:\mu (A)=\mu (B)\neq + \infty ,A,B\in \Cal B
\}.$$
From a functional analytic viewpoint condition $(\Cal L)$ describes a subset
of the kernel of $\mu$ as a linear functional on some subspace of $
L_0$.
Two possibilities for the subspace 
come immediately to mind, $L_1$ and $L_{\infty}.$  If $
\nu$ is a
$\sigma$-finite signed measure satisfying $(\Cal L)$, then taking $A$
to be the
empty set we see that $\nu\ll\mu $ and that if $g=d\nu /d\mu$ then
$\int_{\Omega}fgd\mu =0$ for all of $f\in X$.  If $\nu$ has finite 
variation then $
g\in L_1$ and
the uniqueness question for condition $(\Cal L)$ is equivalent to determining
if $\overline X^{w^{*}}=L^0_{\infty},$ the set of functions h in $
L_{\infty}$ with $\int h d\mu =0.$  If
$d\nu /d\mu\in L_{\infty}$ then uniqueness under $(\Cal L)$ is equivalent 
to $\overline 
X=L^0_1.$  (Here
the overbar denotes norm closure.)

If $\mu$ is $\sigma$-finite but not finite then we need to consider the space of
functions with finite $L_1$ norm on each set of finite measure with the
(locally convex) topology of $L_1$ convergence on sets of finite measure
which we denote by $FL_1.$  The dual of $FL_1$ is
$FL_{\infty}=\{f\in L_{\infty}:\mu (\{f\neq 0\})<\infty \}.$  Note that $
X\subset FL_{\infty}$ and $(\Cal L)$ determines
$\mu$ uniquely if and only if the $\sigma (FL_{\infty},FL_1)$ closure of X in $
FL_{\infty}$ is $FL^0_{\infty}$.
 Let $X_1$ denote the norm closure of X in $L_1$ and let $X_{\infty}$ denote the
$\sigma (FL_{\infty},FL_1)$ closure of X in $FL_{\infty}.$  If $\mu$ is finite $
X_{\infty}=\overline X^{w^{*}}$, i.e., the
$\sigma (L_{\infty},L_1)$ closure.

To see the power of this viewpoint consider the following question
posed in [L].  Is it possible to construct a purely atomic measure
with atoms $(A_i)$ such that $\mu (A_i)\downarrow 0$ and every second atom is 
a bully
but $(\Cal O)$ does uniquely determine $\mu$?  Actually we can show that
neither $A_1$ and $A_2$ can be a bully and have $\mu$ uniquely determined.
Indeed, $A_1$ cannot be a bully since we need only choose $\nu (A_1
)>\mu (A_1)$
and $\nu (A_i)=\mu (A_i)$ for $i>1$.  If $A_1$ is not a bully but $
A_2$ is, define a
map from $L^0_{\infty}$ into $\Bbb R^2$ 
 by $Tf=(\int_{A_1}fd\mu ,\int_{A_2}fd\mu ).$  Then
the range of $T$ is $\Bbb R^2$ but if $f\in X$ then $f=c(1_{A_1}-1_{
A_2})+g$ where supp
$g\subset\cup_{i>2}A_i$.  Thus the range of $T_{|X_{\infty}}$ is 
a subspace properly contained
in $\Bbb R^2$.  This immediately implies that $(\Cal L)$ does not uniquely determine
$\mu$. It will follow from Corollary 1.7 c) below that $(\Cal O)$ cannot
determine $\mu$ uniquely, but we can see this directly as
follows.
Fix $b$, $0<b<1$ and define $\nu 
(A_i)=b\mu (A_i)$ for
$i>2,\nu (A_1)=\mu (A_1),$ and $\nu (A_2)=(1-b)\mu (A_1)+b\mu (A_2
).$  A simple
computation shows that $\nu$ satisfies $(\Cal O).$
The point of the using functional analysis is that it often
reduces questions about $(\Cal L)$ to counting dimensions.

Next we will explore the role of the non-atomic part of $\mu$.
 
\proclaim{Lemma 1.3}If $\mu $ is a finite purely non-atomic positive
measure, then
$X_1=L^0_1$ and $X_{\infty}=L^0_{\infty}$.\endproclaim
\demo{Proof}First note that because $\mu $ is finite $L_{\infty}
\subset L_1$ and thus by
the Hahn-Banach theorem if $X_{\infty}=L^0_{\infty},$ $X_1=
L^0_1.$   Thus we need only
consider the case of $X_{\infty}.$  Observe that $L^0_{\infty}$ is a $
w^{*}$ closed subspace
of $L_{\infty}$ and thus the unit ball of $L^0_{\infty}$ is $
w^{*}$ compact.  It is easy to see 
that an extreme point of the unit ball of $L^0_{\infty}$ is of the form $
1_A-1_B$ where
$\mu (A)=\mu (B)$, $A\cup B=\Omega ,$ and $A\cap B=\emptyset$.  Thus 
all of these extreme points
are in $X_{\infty}$ and it follows from the Krein-Milman theorem that $
X_{\infty}$ is $w^{*}$ 
dense in $L^0_{\infty}$.\qed\enddemo

The argument above used the extreme points which happened to be in
$X$.  If the measure space contains atoms then the extreme points of
the ball of $L^0_{\infty}$ are of the form $1_D-1_E+\lambda 
1_A$ where $D,E\text{ and }A$ are
disjoint, $D\cup E\cup A=\Omega ,A$ is an atom, $|\lambda |\leq 1
,\text{ and }\lambda\mu (A)+\mu (D)=\mu (E)$.
Whether or not such things are in $X_{\infty}$ is a matter of the
combinatorics of the measures of the atoms, however, we can still
say something about the non-atomic part.
\proclaim{Proposition 1.4}If $\mu$ is a $\sigma$-finite positive
measure on $(\Omega,\Cal B)$, $\mu =\mu_c+\mu_a$, where $\mu_c$ has 
no atoms and $
\mu_a$ is purely
atomic, and $\nu$ is a $\sigma$-finite signed measure which satisfies $
(\Cal L)$ then
$\dsize\frac{d\nu}{d\mu_c}$ is constant.\endproclaim
\demo{Proof}The hypothesis implies that $\nu\ll\mu $ and thus if $
A\in \Cal B$
such that $\mu (A)$ and $\nu (A)$ are finite, $g1_A\in L_1,$ where $
\dsize g=\frac{d\nu}{d\mu}.$  Because
both $\mu$ and $\nu$ are $\sigma$-finite we can find a sequence of sets $
D_n\in \Cal B$
such that $\mu (D_n)$ is finite and non-zero, $\cup_nD_n$ contains the
complement of the union of the atoms, and $g_{|D_n}\in L_{\infty}(
\mu_{|D_n})$ for all n.
Now by Lemma 1.3, $X\cap L^0_1({\mu_c}_{|D_1\cup D_j})$ is dense in $
L^0_1({\mu_c}_{|D_1\cup D_j})$ for
each j.  Thus $g$ is constant on $D_1\cup D_j$ for all j and therefore $
g$ is
constant on $\cup D_n$, establishing the result.\qed\enddemo

Proposition 1.4 shows that the difficulty really is in the atomic part of 
$\mu$.  Next we show that uniqueness can be established provided it holds 
for a rich enough family of restrictions of $\mu$.
 
\proclaim{Proposition 1.5}Suppose that $\mu$ is a $\sigma$-finite measure
and for every pair of disjoint sets $D$ and $E$ of finite measure there
is a measurable set $F$ such that $X_{\infty}\cap FL_{\infty}({\mu}_{|F})
=FL^0_{\infty}({\mu}_{|F}),$
respectively, $X_1\cap FL_1({\mu}_{|F})=FL^0_1({\mu}_{|F}),$
and there exist $D^{\prime}$ and $E^{\prime}$
contained in F such that $\mu (D)=\mu (D^{\prime})$ and $\mu (E)=\mu 
(E^{\prime})$.  Then if
$g\in FL_1$, respectively, $g\in L_{\infty}(\mu ),$ and $g|_X=0,$ then 
$g$ is constant.\endproclaim
 
\demo{Proof}Because $\text{sp}\{\mu (D)1_E-\mu (E)1_D:D,E\text{ disjoint 
and finite measure}
\}$
is dense in $FL^0_{\infty}$ and in $FL^0_1$, 
in either case it is sufficient to show
that $\mu (E)\int_Dg d\mu =\mu (D)\int_Eg d\mu$ 
for all disjoint sets and $
D$ and $E$ of
finite measure. 
So given $D$ and $E$ let $F$, $D^{\prime}$, and $
E^{\prime}$
satisfy the hypothesis of the proposition. 
Because $\mu (D')1_{E'}-\mu (E')1_{D'},$ $1_{D'}-1_D,$ and $1_{E'}-1_E$ are
in $X$, 
 $\mu (E')\int_{D'}g d\mu =\mu (D')\int_{E'}g d\mu$,
$\int_{D'}g d\mu=\int_{D}g d\mu$, and $\int_{E'}g
d\mu=\int_{E}g d\mu$. Hence $$\mu (E)\int_Dg d\mu =\mu
(E')\int_{D'}g d\mu =\mu (D')\int_{E'}g d\mu=\mu (D)\int_Eg
d\mu.\qed$$  \enddemo

In the papers [L] and [N] the following auxiliary function plays an
important role.  Let $W=\text{range of }\mu$ and suppose that $\nu$ satisfies $
(\Cal L)$.
For each $w\in W$ define $f(w)=\nu (A)\text{ if }\mu (A)=w$. The
absolute continuity of $\nu$ with respect to $\mu$ implies that
$f(0) =0.$   The lemma below
summarizes the important properties of $f$.  Versions of a)-c) are in
[L] and d) and e) are in [N].  Because the arguments are short we
include them for the sake of completeness.
 
\proclaim{Lemma 1.6}Suppose that $\mu$ is a  $\sigma$-finite
positive measure and $\nu$ is a $\sigma$-finite signed
measure on $(\Omega ,\Cal B)$. If $\nu$ satisfies $(\Cal L),$  then with $
f$ as above
\roster
\item "a)"$f$ is well defined and finite for all $w$ finite in W.
\item "b)"$(\Cal O)$ is equivalent to the assertion that $f$ is non-decreasing.
\item "c)" If there exists sets $A\in \Cal B$ of arbitrarily
small but positive measure, and  $f^{\prime}(w)$ exists at some $w\in
W$, then
$\dsize\lim_{\mu (A)\to 0}\frac{\nu (A)}{\mu (A)}$ exists and equals $f^{\prime}(w).$ 
$(f^{\prime}$ is 
defined only at limit points of W.)
In particular, if $(\Cal O)$ holds and W has positive measure then
$\dsize\lim_{\mu (A)\to 0}\frac{\nu (A)}{\mu (A)}$ exists.
\endroster
If range $\mu$ is a finite interval, we have in addition
\roster
\item "d)"$f$ is continuous.
\item "e)"If $\nu\geq 0$, then $f$ is nondecreasing. Consequently, $(\Cal O)$ is
satisfied  and $f'(w)$ exists  and is constant
a.e. $\mu$.
\endroster
\endproclaim
\demo{Proof}a) and b) are obvious. 
Let $w$
be a point in the range of $\mu$ where the derivative exists and suppose
that $\mu (B)=w$.  If $A_n\subset B,\mu (A_n)\neq 0$, and $\mu (
A_n)\rightarrow 0$ then
$$f^{\prime}(w)=\lim\frac {\nu (B\setminus A_n)-\nu (B)}{\mu (B\setminus 
A_n)-\mu (B)}=\lim\frac {-\nu (A_n)}{-\mu (A_n)}.$$
If $A_n\cap B=\emptyset ,\mu (A_n)\neq 0$ and $\mu (A_n)\rightarrow 
0,$ then
$$f^{\prime}(w)=\lim\frac {\nu (B\cup A_n)-\nu (B)}{\mu ( B\cup A_n
)-\mu (B)}=\lim\frac {\nu (A_n)}{\mu (A_n)}.$$
Finally if $(A_n)$ is any sequence of sets with $\mu (A_n)\rightarrow 
0$,
$$\nu (A_n)-f^{\prime}(w)\mu (A_n)=\nu (A_n\cap B)-f^{\prime}(w)\mu 
(A_n\cap B)+\nu (A_n\cap B^c)-f^{\prime}(w)\mu (A_nB^c)$$
$$=o(\mu (A_n\cap B))+o(\mu (A_n\cap B^c))=o(\mu (A_n)).$$
If $W$ has positive measure there must be a sequence of sets
$(A_n)$ with $\mu (A_n)\rightarrow 0$.  Indeed, if not, $\mu$ 
would be purely atomic with the measures of the atoms bounded
away from 0 and thus
the range of $\mu$ would be a countable set. If $(\Cal
O)$ is satisfied then $f'$ exists at almost every point of $W$.

For d) we will show that continuity at $\mu (B)=x$ follows from the
absolute continuity of $\nu$ with respect to $\mu$.  First if $\mu_
c\neq 0,$ by
Proposition 1.4 there is a constant $\kappa >0$ such that $\nu =
\kappa\mu_c$ or $ -\kappa\mu_c$ on $\Omega\setminus\cup A_i,$
where $(A_i)$ is the sequence of atoms of $\mu$.  Given $\varepsilon 
>0$ find
$\delta ,\ \varepsilon /(3\kappa )>\delta >0$ such that $|\nu| (A)<\varepsilon 
/6$ if $\mu (A)<\delta .$  Choose n such that
$\sum_{i\geq n}\mu (A_i)<\delta /2.$ Let $\{r_j\}$ be the range of $
\mu |_{\cup_{i<n}A_i}$ and let $K$ be
the range of $\mu |_{C\cup\cup_{i\geq n}A_i}$.  Because $K$ 
is compact there is a number
$\delta^{\prime}>0$ such that if $x\notin r_i+K,(x-\delta^{\prime}
,x+\delta^{\prime})\cap r_i+K=\emptyset$.  Let
$\rho =\text{min}\{\delta ,\varepsilon /(3\kappa ),\delta^{\prime}
\}$.  If $y\in W$ and $|x-y|<\rho$,  then for some $i$, $x$ and
$y$ are in $r_i+K.$  Therefore
$$|f(x)-f(y)|=|f(r_i)+\nu (\cup_{i\in I}A_i)+\nu (D)-(f(r_i)+\nu (\cup_{
i\in J}A_j)+\nu (E))|$$
where $I\cup J\subset \{i\geq n\}$, $D\cup E\subset C,$ $x=r_i+\mu
(\cup_{i\in I}A_i\cup D)$, and
$y=r_i+\mu (\cup_{i\in J}A_i\cup E).$  Because $\mu (\cup_{i\in
I\cup J}A_i)<\delta /2$ and
$|\nu (D)-\nu (E)|=\kappa |\mu (D)-\mu (E)|$,
$$|f(x)-f(y)|<\varepsilon /3+\kappa |\mu (D)-\mu (E)|\leq \varepsilon
/3+\kappa [|x-y|+\delta ]<\varepsilon .$$

To see that $f$ is nondecreasing when $\mu$, $\nu\geq 0$, suppose that this 
is false.
Then $f$ must have a local minimum at some point $x$, $0<x<\mu (\Omega
)$, and
we may assume that $f(y)>f(x)$ for all y in some interval $(x-\delta
,x),$
$\delta >0$.  Because the range of a finite measure is compact, either
$x=\mu (\cup_{i\in I}A_i)+\mu (D)$ with $D\subset C$ and $\mu (D)>
0$, or $x=\sup \{\sum_{i\in I}\mu (A_i):$
$\sum_{i\in I}\mu (A_i)<x\}$ for some infinite set $I\subset \Bbb N$.  
In either case $
f(x)=f(y_n)+\nu (B_n)$
for some sequence $(y_n)$ increasing to $x$.  Since $\nu (B_n)\geq
0$, $f(x)$ cannot
be a strict local (left) minimum.
\qed\enddemo

The results in [L] and [N] were obtained by showing that in the case in
which $W$ is an interval $f^{\prime}$ is constant and there is
no singular part to $\nu$ and thus $f(x)=x$.  The
function $f$ is unsatisfactory for generalizing the results.  We will
usually use the Radon-Nikodym derivative of $\nu$ with respect to $
\mu$ in
our arguments.  The Radon-Nikodym derivative is more indicative of
the situation because it retains the information about the underlying
sets whereas $f$ mostly reflects properties of the range.
 
\proclaim{Corollary 1.7}
\roster
\item "a)"If $\mu$ is a $\sigma$-finite measure with range of positive
Lebesgue measure such that for every $\varepsilon >0,\mu (\cup \{A
:A$ is an atom and
$\mu (A)>\varepsilon \})<\infty$ and $\nu$ is a $\sigma$-finite signed measure satisfying $
(\Cal O)$, then
$\dsize\frac{d\nu}{d\mu}
\in L_{\infty}(\mu )$.
\item "b)"If $\mu$ is a positive measure with range a finite interval and $
\nu$
is a positive measure satisfying $(\Cal L)$, then $\dsize\frac{d\nu}{d\mu}
\in L_{\infty}
(\mu )$.
\item "c)" If $\mu$ is  a positive $\sigma$-finite measure 
, $\nu$
satisfies $(\Cal L)$ and $\dsize\frac{d\nu}{d\mu} \in L_\infty(\mu)$
and is not constant then there is a $\sigma$-finite 
positive measure 
$\nu'$ satisfying $
(\Cal L)$ and $\dsize\frac{d\nu}{d\mu}$ is not constant. If, in
addition, the range of $\mu$ is a finite interval, then $\nu'$
satisfies $(\Cal O)$.

\endroster
\endproclaim
 
\demo{Proof}By Proposition 1.4, $\dsize\frac{d\nu}{d\mu}
$ is constant on $C$, the complement
of the atoms, and by c) of Lemma 1.6, $\dsize\lim_{\mu (A)\to 0}\frac{\nu
(A)}{\mu (A)}
$ exists.  Thus
$\dsize\frac{\nu (A)}{\mu (A)}$ is bounded for all atoms with $\mu (A)$ 
sufficiently small.  The
hypothesis implies that there are only finitely many other atoms,
and thus $\dsize\frac{d\nu}{d\mu}
$ is bounded.  Part b) follows from e) of Lemma 1.6
and a).
 If $\mu\geq 0$ and $\dsize\frac{d\nu}{d\mu}
\in L_{\infty}(
\mu )$ is not constant and $(\Cal L)$
holds then there is a nonzero constant $\gamma$ such that
$\nu^{\prime}
=\mu +\gamma\nu\geq 0.$
Clearly $\nu^{\prime}$ also satisfies $(\Cal L)$. If the range of
$\mu$ is an interval, by Lemma 1.6 e)
 $\nu'$ satisfies $(\Cal O)$.
\qed\enddemo
 
Note that Corollary 1.7 c) says that if the range of $\mu$ is
a finite interval, then the uniqueness question for
signed measures with bounded Radon-Nikodym derivative is the same
for $(\Cal L)$ and $(\Cal O)$.

\heading  2. Finite Signed Measures \endheading
Next we will prove some technical results which are useful for
establishing uniqueness for the case of signed measures.
 
\proclaim{Lemma 2.1}Suppose that $\mu$ is a finite positive measure with
atoms $\{A_n\}$ arranged so that $\mu (A_n)\geq \mu (A_{n+1})$ and let $
C$ be the
complement of the union of the atoms.  If $Z$ is a subspace of $L_{
\infty}$
which contains $\{1_{A_n}-\mu (A_n)/(\mu (C\cup\cup_{j>n}A_j)
1_{C\cup\cup_{j>n}A_j}:n=1,2,...\}$ \newline
and $L^0_{\infty}(C,\mu)$, 
then $\overline
Z^{w^{*}}\supset L^0_{\infty}$.
\endproclaim
 
\demo{Proof}It is sufficient to show that
$$\frak D=\{\mu (A_1)1_{A_k}-\mu (A_k)1_{A_1}:k=1,2,...\}\cup \{\mu
(C)1_{A_1}-\mu (A_1)1_C\}\subset\overline Z^{w^{*}}.$$
Indeed if $g\in L_1$ and $z(g)=0$ for all $z\in \frak D$, then $\mu(A_k)
g(A_1)=\mu(A_1)g(A_k)$, for all k and
$g(A_1)\mu(C)=\int_Cgd\mu$.  Because $L^0_{\infty}
(C,\mu)\subset Z$, $g$ is a constant $\gamma$ 
on $C$.  It follows that $g$ is $\gamma$ on $
\Omega$.

For each $n$ let $B_n=C\cup\cup_{k>n}A_n$.  In this notation we have from the
hypothesis that $1_{A_n}-\mu (A_n)/\mu(B_n)1_{B_n}\in Z$ for all $
n$.  Now observe that since
$B_1=A_2\cup B_2$,
$$\multline 1_{A_1}-\mu (A_1)/\mu(B_1)1_{B_1}+
\mu (A_1)/\mu(B_1)(1_{A_2}-\mu (A_2)/\mu(B_2)1_{B_2}
)\\
=1_{A_1}-\mu (A_1)/\mu(B_2)1_{B_2}\in Z .
\endmultline
$$
An easy induction argument shows that $1_{A_1}-\mu (A_1)/\mu(B_
k)1_{B_k}\in Z$ for all k.
Therefore
$$\multline 1_{A_1}-\mu (A_1)/\mu (A_n)1_{A_n}=1_{A_1}-\mu (A_1)/
\mu (B_{n-1})1_{B_{n-1}}\\
+(\mu (A_1)/\mu(B_{n-1})-\mu (A_1)/\mu (A_n))(1_{A_
n}-\mu (A_n)/\mu(B_n)1_{B_n})\in Z\endmultline$$
for all $n$.  If $\mu (C)=0$ or there are only finitely many atoms, then
our initial observation shows that $\overline Z^{w^{*}}\supset L^{
0}_{\infty}$.  If $\mu(C) \neq 0$ and there
are infinitely many atoms, then note that $||1_{A_1}-\mu (A_1)/\mu
(B_n)1_{B_n}||_{\infty}\leq$
$\max \{1,\mu (A_1)/\mu(C)\}$ for all $n$ and that $1_{A_1}-\mu
(A_1)/\mu(B_n)1_{B_n}\rightarrow 1_{A_1}-\mu (A_1)/\mu
(C)1_C$
in the $w^{*}$ topology.  Again the conclusion follows from our initial
observation.\qed\enddemo
 
\proclaim{Theorem 2.2}Suppose that $\mu$ is a finite positive measure with
atoms $\{A_i\},$ $\mu(A_i)\geq \mu(A_{i+1}),$ and there is a constant $
K$ such that
for each $i,$ $m\in \Bbb N$ there is an $h\in X_1$ with $\text{supp
}
h-1_{
A_i}\subset\Omega\setminus\cup_{j\leq m}A_j$ and
$||h||_1\leq K\mu(A_i).$  If $F\in X^{\perp}_1\cap L_{\infty}$ 
and $\lim_{
n\to\infty}$ $F(A_n)$ exists, then $F$
is a constant.\endproclaim
 
\demo{Proof}We will show the hypothesis of Lemma 2.1 is satisfied
where $Z=\{z\in L_{\infty}:\int zFd\mu=0\}$.

Let $C$ be the complement of the union of the atoms.  Fix $n\in \Bbb N$.  Let
$g=1_{A_n}-\mu (A_n)/(\mu(C\cup\cup_{j>n}A_j$ $)1_{C\cup
\cup_{j>n}A_j}$.  For each integer $m>n$
let
$$g_m=h_{mn}-\mu (A_n)/\mu(C\cup\cup_{j>n}A_j)\sum^m_{s=n+1}
h_{ms}$$
where $h_{ms}\in X_{\infty}$ satisfies $||h_{ms}||_1\leq K\mu(
A_s)$, $h_{ms|\cup_{j<s}A_j}=0$, $h_{ms|A_s}$
$=1_{A_s},$ and $h_{ms|\cup_{s<j\leq m}A_j}=0$, for $n\leq
s\leq m$.  Therefore
$$||g_m||_1\leq K\left[\mu(A_n)+\mu(A_n)/\mu(C\cup\cup_{
j>n}A_j)\sum^m_{s=n+1}\mu(A_s)\right]\leq 2K\mu(A_n).$$

Suppose that $F\in L_{\infty}\cap X^{\perp}_1$ and $\lim 
F(A_n)=\rho$.  If $\mu(C)>0$,  by
Proposition 1.4, $L^0_1(C,\mu)$ is contained in $X_1$ and thus $
F$ is constant on $C$.
Also note that $F(c)=\rho$ for all $c\in C,$ because $
F$ is constant on $C$
and if $0<\mu (A_i)\leq \mu(C),$ 
$1_{C^{^{\prime}}}-1_{A_i}\in X_1$ where $C^{^{\prime}}
\subset C$ such that
$\mu (C^{^{\prime}})=\mu (A_i).$  Hence
$$\align\left|\int (g-g_m)Fd\mu\right|&=\left|\int_{\cup_{j>m}A_
j\cup C}(g-g_m)Fd\mu\right|\\
&\leq ||(F-\rho )|_{(\cup_{j>m}A_j\cup C)}||_{\infty}||g-g_
m||_1+|\int (g-g_m)\rho d\mu |.\\
\endalign$$
 
The second term is zero because $g$ and $g_m\in L^0_{\infty}$ 
and the first term
clearly converges to zero. Since $g_m \in X_1$, $\int g_m F
d\mu=0$, for each $m$. Thus $\int gFd\mu =0$ and therefore $
\{F\}^{\perp}\supset L^0_{\infty}$
and $F$ is constant as claimed.\qed\enddemo

As an application of Theorem 2.2 we will prove Nymann's Theorem.  We
will need the following lemma whose proof we leave to the reader.
\proclaim{Lemma 2.3}If f and $g\in L_1(\mu )$ for some measure $\mu$ and there
is a set A such that $f_{|A}=\gamma g_{|A}$ for some $\gamma\in [-
1,0]$ and
$||f_{|A}||_1\geq ||f-f_{|A}||_1,\text{ then }||f+g||_1\leq ||g|
|_1$.\endproclaim
 
\proclaim{Theorem 2.4 }If $\mu$ is a positive measure with range
equal to a finite interval, $\nu$ is a measure satisfying $
(\Cal L)$, and $\nu$ is positive or $\nu$ is a signed measure with
$\dsize\frac{d\nu}{d\mu}\in L_\infty(\mu)$,
then $\nu$ is a constant times $\mu$.\endproclaim
 
\demo{Proof} Assume first that $\nu$ is positive. We will verify
the hypotheses of Theorem 2.2 for
$F=d\nu /d\mu$.  By Lemma 1.6, $\lim_{
\mu (A)\to 0}\nu (A)/\mu (A)=c$ 
so $\lim_{n\rightarrow\infty} F(t_n)=c$,  if $t_n\in A_n$ and $\mu (A_n)\to 0$.  
Also $F$ is a constant on $
C$ by 
Proposition 1.4 and therefore $F\in L_{\infty}(\mu )$.  
Thus we need only produce 
the functions $h$ in $X_1.$

Because the range of $\mu$ is an interval for each $n$ there exists a
subset $M_n$ of $\{n+1,n+2,...\}$ and $C_n\subset C$ such that
$\mu (A_n)=\mu (C_n\cup\cup_{j\in M_n}A_j)$,  that is, $h_n=1_{A_n}
-1_{C_n\cup\cup_{j\in M_n}A_j}\in X_1$
$\text{and }||h_{n|A_n}||=||h_n-h_{n|A_n}||$. (See [L, Proposition
1].)
\example{CLAIM}  If $k$ and $m$ are integers, $k<m,$ and $g\in X_1$, then there is an 
$h\in X_1$ such that $h_{|\cup_{n\leq k}A_n}=g_{|\cup_{n\leq k}A_
n},\text{ }h_{|\cup_{k<n\leq m}A_n}=0,$ and  $||g||_{L_{1(\mu )}}\geq 
||h||_{L_{1(\mu )}}$.
\endexample
 
Indeed, suppose that $g|_{A_{k+1}}=a1_{A_{k+1}}.$  Let $f=-ah_n$ in Lemma 2.3.
Then if $g_1=g+f,g_1|_{\cup_{n\leq k}A_n}=g_{|\cup_{n\leq k}A_n}
,g|_{A_{k+1}}=0,$ and $||g_1||\leq ||g||$.  Induction
finishes the proof of the claim.  Applying the claim to $g=h_n$ shows
that the hypothesis of Theorem 2.2 is satisfied and thus $d\nu /d\mu$ is
constant.

For $\nu$ signed, we apply the previous argument to $\nu'$ from
Corollary 1.7 c). It follows that $\nu'$ and hence $\nu$ is a
multiple of $\mu$. \qed\enddemo

\remark{Remark 2.5} This theorem includes Nymann's result. The
proof given by Nymann uses the result of Malitz [M] which in
turn uses [L]. Thus one virtue of the proof we have given is that
it is direct.  More important our argument shows that there are
two main ingredients to the proof. The ability to solve equations
and the ``continuity at $\emptyset$'' of the Radon-Nikodym
derivative.

This second ingredient is also the main obstacle to be overcome
if the requirements on signed measures are to be relaxed.
 It would
be nice to find a replacement for the argument used in the proof of
Corollary 1.7 so that signed measures without bounded
Radon-Nikodym derivative could also be treated.  In
particular we do not know whether a finite signed measure with
range an interval is determined uniquely by $(\Cal L)$.
\endremark
 
\heading {\bf 3.  $\sigma$-Finite Measures}\endheading

In the paper [N] only finite positive measures were considered.  Here we are
allowing $\sigma$-finite measures.  At first glance it may appear that the
$\sigma$-finite case should reduce to the finite case.  However we wish
to observe that given $t$ $\in$ range $\mu$ there may be no set $C$ of finite
measure such that range $\mu_{|C}$ is an interval containing $t$.  Consider the
following example.

\example{Example}  Let $\Omega$ be [0,1] union a sequence of disjoint sets $
\{A_n\}$
(also disjoint from [0,1]) and let $\mu$ restricted to [0,1] be Lebesque
measure and let the $A_n$'$s$ be atoms of measure $1+2^{-n}$.  Note that the
range of $\mu$ is $\Bbb R^{+}$ but for every $C$ of finite measure
greater than one
the range of $\mu_{|C}$ misses some interval $(1,1+2^{-n}),$ with
$n\in \Bbb N$.\endexample
 
Thus it is not possible to reduce the problem to the finite case
directly.  The next few lemmas will show that the example given
above is the prototype for this difficulty.
 
\proclaim{Lemma 3.1}Let $\mu$ be a purely atomic measure with range equal
to an interval such that for every $\varepsilon >0,\mu (\cup \{A_n
:A_n$ is an atom and
$\mu (A_n)>\varepsilon \})<\infty$.  Then for every $t< \mu (\Omega
)$ there is a subset $C$ of $\Omega$
such that $\mu (C)\geq t$ and the range of $\mu_{|C}$ is a finite 
interval.\endproclaim
 
\demo{Proof}There is nothing to prove if $\mu (\Omega )<\infty ,$ so assume that
$\mu (\Omega )=\infty$.  The requirement that the range is equal to a finite
interval is
equivalent to the assumption that there are no bullies.  Let $(a_n
)$ be
the masses at the atoms arranged in decreasing order.  Assume that
$a_1=1$.  We will inductively construct a subsequence of the $a_n$'s
without bullies so that the sum is finite.

As a convenience we will say that if $E$ and $F$ are subsets of $\Bbb N$ and
every element of $E$ is less than every element of $F$ then $F>E$.  Let
$E_0=\{1,2,...N\}$ where $N$ satisfies $\sum_{i\leq N} a_i> t$.  Choose
$E_1>E_0$ such that $\frac 12\leq \sum_{n\in
E_1}a_n\leq\frac 23$ and $a_n\leq\frac 12$
for all $n\in E_1$.  This is possible because $\sum_{n\geq k}a_
n=\infty$ for all $k$ and $(a_k)$
decreases to 0.  Next choose $E_2>E_1$ such that $\frac 13\geq \sum_{
n\in E_2}a_n\geq\frac 14$
and $a_n\leq\frac 14$ for all $n\in E_2$.  Thus we can continue in this way to
construct a sequence of finite subsets of $\Bbb N$,
$(E_j),$ such that for all $j \geq 1$
\roster
\item "a)"if $n\in E_j,a_n\leq 2^{-j}$
\item "b)"$\displaystyle{2^{-j+2}/3\geq \sum_{n\in E_j}a_n\geq 2^{
-j}}$
\item "c)"$E_{j+1}>E_j$
\endroster
 
It is easy to see that $C=\cup \{A_n:n\in E_j$, some $j\}$ is the required set.
If $n\in E_j$, then $a_n=\mu (A_n)\leq 2^{-j}$ and $\sum_{k>j}\sum_{
n\in E_k}a_n\geq \sum_{k>j}2^{-k}=2^{-j}$.
Thus $A_n$ is not a bully.
\qed\enddemo
 
\proclaim{Lemma 3.2}If $\mu$ is a positive
$\sigma$-finite measure with range equal to
an interval, then there is a possibly infinite number $\beta$ and a
subset $B\text{ of }\Omega$ such that
\roster
\item "a)"The range of $\mu_{|B}$ is the interval $[0,\beta]$ and for every $
t<\beta$
there is a set $B^{\prime}\subset B$ with the 
range of $\mu_{|B^{\prime}}$ an interval of finite
length at least t.
\item "b)"either
\itemitem{i)} $B=\Omega$
\item "\hphantom{b)}"or
\itemitem{ii)} $\mu_{|\Omega\setminus B}$ is purely atomic, every atom in $
\Omega\setminus B$ has measure
strictly greater than $\beta$ and $\beta$ is a limit point of the range of $
\mu_{|\Omega\setminus B}$.
\endroster
\endproclaim
 
\demo{Proof}The result is trivial if $\mu (\Omega )<\infty .$  If $
\mu (\Omega )=\infty$ there are
two cases to consider.
\roster
\item "1)"there is a sequence of atoms $(A_i)$ such that $\lim \mu
(A_i)=0$ and 
$\sum\mu (A_i)=\infty$
\item "2)"there is no such sequence.
\endroster

If Case 1 occurs, let $A=\cup A_i$, then range of $\mu_{|A}=\Bbb R^{
+}$ and by Lemma 3.1
we can take $B=\Omega$ and $\beta$ $=\infty$.

Now we will treat Case 2.  As usual let $C$ be the complement of the atoms
and let $\kappa =\mu (C)$, which we assume is finite.  (If $\kappa
=\infty$, we can
again let $B=\Omega$.)  Let $\bold B=\{G:G\subset\Omega\setminus
C$, $\mu(G)<\infty$ and range $\mu_{|G}$ omits no
intervals of length greater than $\kappa \}$.  Note that if $B=\cup_{
i\in M}G_i,$ $G_i\in \bold B$,
for each $i\in M$, and $\mu (B)<\infty$, then range $\mu_{|B}$ 
omits no intervals of length
greater than $\kappa$.  Indeed, if $H\subset B$ and $s = \sup \{\mu (D)
:D\subset B,$ $\mu (D)<\mu (H)\}<\mu(H)-\kappa$,
then every non-null $E\subset H$ has measure greater
than $\kappa$.  Let E be an atom of
minimal measure contained in $H$, then there exists $G_i\supset E$ with range
$\mu_{|G_i}$ missing no interval of length greater than $\kappa$.  But for some
$t\in [\mu (E)-\kappa ,\mu (E))$ there is a subset F of $G_i$ so that $
\mu (F)=t$.  Clearly
$F\cap H=\emptyset$ so $\{\mu ((H$ $\setminus$ $E)\cup F):F\subset
G\}$ contains a point in $(s,\mu (H))$.
Hence $\mu_{|B}$ omits no intervals of length greater than $\kappa$.

It follows that if $B_1=\cup \bold B\cup C$, then if $t<\mu (\cup
\bold B)+\kappa$ there is a subset
$B^{\prime}$ of $\cup \bold B$ such that range $\mu_{|\text{$B^{\prime}$ $
\cup C$}}$ is an interval containing [0,t] and
$\mu (\text{$B^{\prime}\text{$\cup C$ $)<\infty$}$}$.  If $\mu (B_
1)=\infty$, then let $B=\Omega$.  Otherwise let $B=B_1$ and
$\beta =\mu (\text{$\cup \bold B)+$}\kappa$.  Also note that $\beta
>0$.  Indeed if $\beta =0,1)$ Case 1 does not occur and there is
some $\varepsilon >0$ so that $\mu (F)<\infty$ 
where $F=\cup \{A:A$ is an atom and $
\mu (A)<\varepsilon \}$.  Because
the range of $\mu$ is an interval, $F\neq\emptyset$.  The range of $
\mu_{|F}$ cannot be an
interval because $\beta =0$, so there must be a bully in $F$.  However
this implies that the range of $\mu$ is not an interval.

Clearly the range of $\mu_{|B}=[0,\beta]$.
Observe that if $A$ is an atom with $\mu (A)\leq \beta <\infty$ then $
A\cup B\in \bold B$ and
hence $A\subset B$.  Therefore if $\beta <\infty$ and the range of $
\mu$ is $[0,\infty ]$, every atom
$A\subset\Omega$ $\setminus$ $B$ has measure greater than $\beta$ 
and there must be an infinite sequence
of atoms $(A_i)$ such that $\mu (A_i)$ decreases to $\beta$.
\qed\enddemo

We are now in a position to extend the result in [N] to all
positive	
$\sigma$-finite measures with range an interval.
 
\proclaim{Proposition 3.3}Suppose that $\mu$ is a  positive $\sigma$-finite 
measure  with range an interval and that $\nu$ is a
$\sigma$-finite
measure satisfying
$(\Cal L)$. If $\nu$ is positive or if there exist $\epsilon >0$
and $M<\infty$ such that $\dsize\frac{|\nu(A)|}{\mu(A)}<M$ if
$\mu(A)<\epsilon $, then $\nu$ is a multiple of $
\mu$.\endproclaim
 
\demo{Proof}By Theorem 2.4, either condition implies that if
$B\in \Cal B$, $\mu(B)<\infty$ and range $\mu_{|B}$ is an
interval, then $\dsize\frac{d\nu}{d\mu}$ is constant on $B$.
Hence we assume that
$\mu (\Omega )=\infty$.  From Lemma 3.2 it follows that if $\beta
=\infty$ then for any
pair of disjoint sets of finite measure there is a subset $B^{^{\prime}} 
\subset B$ such
$\mu_{|B^{^{\prime}}}$ has range a finite interval containing $\mu
(D)$ and $\mu (E)$.  By
Proposition 1.5 and the finite case this implies that $d\nu /d\mu$ is
constant.

Now assume that $\beta <\infty$.  In view of the Proposition 1.5 it is
sufficient to show that for $r>s>0$ there is a set $F$ such that
$\overline {X\cap L_{\text{$\infty$}}(\mu_{|F})}^{w^{*}}=L^0_{\text{$
\infty$}}(\mu_{|F})$, the range of $\mu_{|F}$ contains $\{s,r\}$ and $
\mu (F)<\infty$.
If $\beta\geq r$ we can take $F=B$ and proceed as above.  If $r>\beta$, we
will add some atoms to $B$ to get the required set $F$.  From Lemma
3.2 we know that there is an infinite sequence of atoms $(A_i)$ such
that $\mu (A_i)$ decreases to $\beta$.  It follows that there 
is a finite set of the
$A_i$'$s$, $B_1,$ $B_2$, ... $B_n$, such that $\mu (B_i)>\text{$\mu
(B_{i+1})$}$ for all $i$,
$\mu (B_1)+\beta >\text{$\mu (B_2)+\text{$\mu (B_3)$ $>$ $\mu (B_1
)$}$}$, and both $r$ and $s$ are in the range of
$\mu_{|F}$, where $F=B\cup\cup_{i\leq n}B_i$.  We will show that $
X_{\infty}$ is $w^{*}$ dense in $L^0_{\text{$\infty$}}(\mu_{
|F})$.
Because the range of $\mu_{|B}$ is a finite interval, $X_{\infty}\cap\text{$
L^0_{\text{$\infty$}}(\mu_{|B})$}$ is $w^{*}$
dense in $L^0_{\text{$\infty$}}(\mu_{|B})$.  Also 
$\text{dim } L^0_{\text{$
\infty$}}(\mu_{|F})$/$L^0_{\text{$\infty$}}(\mu_{|B})=n$, 
so it is sufficient to find n
linearly independent elements of this quotient which arise from X.
For each $i<n$ there is a subset $C_i$ of $B$ with positive measure such
that $1_{B_i}-1_{B_n}-1_{C_i}\in X$, and there is a set $C_n\subset
B$ such that $1_{B_1}+1_{C_n}-1_{B_2}-1_{B_3}\in X.$
It is easy to see that these are in $n$ linearly independent cosets of
$\overline {X\cap L_{\infty}(\mu_{|F})}^{w^{*}}/L^0_{\infty}
(\mu_{|B})$.
\qed\enddemo
 
\remark{Remark 3.4} Recently Khamsi and Nymann [KN] have shown that
if $\mu$ is a finite measure with range equal to an interval then any
$\sigma$-finite measure $\nu$ such that if $\mu (A)=\mu (B)$ then $
\nu (A)=\nu (B)$, i.e., $(\Cal L)$
without the finiteness requirement is satisfied, then $(\Cal L)$ holds.  It
follows that Proposition 3.3 holds with the weakened version of $(
\Cal L)$
as well.\endremark
 
\proclaim{Theorem 3.5} Suppose that $\mu$ is a  $\sigma$-finite
signed measure  with range an interval and that $\nu$ is a
$\sigma$-finite
measure satisfying
$(\Cal L)$. If there exist $\epsilon >0$
and $M<\infty$ such that $\dsize\left|\frac{\nu(A)}{\mu(A)}\right|<M$ if
$\mu(A)<\epsilon $ and $A$ is an atom, then $\nu$ is a multiple of $
\mu$.\endproclaim
\demo{Proof} By Lemmas 1.1 and 1.2, $|\mu|$ and $\nu'$ also satisfy the
hypothesis. By Proposition 3.3, $\nu'$ is a multiple of $|\mu|$
and hence $\nu$ is a multiple of $\mu.$ \qed\enddemo

\proclaim{Corollary 3.6}  Suppose that $\mu$ and $\nu$ are 
$\sigma$-finite
signed measures satisfying $(\Cal O)$. If the range of $\mu$ is an interval,
then $\nu$ is a  multiple of $\mu$.
\endproclaim
\demo{Proof} Proceeding as in the previous proof, $|\mu|$ and
$\nu'$ satisfy $(\Cal O)$ and by Lemma 1.6 b), $\nu'$ is positive.
Therefore by  Proposition 3.3, $\nu'$ is a multiple of $|\mu|$
and hence $\nu$ is a multiple of $\mu.$ \qed\enddemo

\heading {\bf 4.  Examples.}\endheading

In order to construct examples we will use the finite dimensional
case to control the combinatorics in the same spirit as the examples
constructed by Leth.

Let us now consider the finite dimensional case, i.e., suppose that
there are finitely many atoms $(A_i)^n_{i=1}$ with $\mu (A_i)=a_i$ and assume
that $\Omega =\cup A_i$.  The statement of the problem can be rephrased in
terms of dimension and thus the uniqueness question for $(\Cal L)$ reduces
to linear algebra.  To make the statement more succinct let us
introduce the following notation.  Let
$\Cal R=\{\nu\in \{-1,0,1\}^n:\langle\nu ,(a_i)\rangle =0\}$.
 
\proclaim{Proposition 4.1}Condition $(\Cal L)$ uniquely determines $\mu$ if and
only if dim sp $\Cal R =n-1$.\endproclaim

It would be nice if there were some easy characterization of the
finite sequences $(a_i)$ so that dim sp $\Cal R=n-1$, but as we will see
below there is no obvious description.  One easy observation is that
the $a_i$'s must all belong to the same rational equivalence class.  (This
same observation also applies in the case where $\mu$ is purely atomic
$\sigma$-finite but the measures of the atoms are bounded away from 0.)
Note also that uniqueness for $(\Cal O)$ implies
uniqueness for $(\Cal L)$ in the
case of finitely many atoms because the Radon-Nikodym derivative is
automatically bounded.

For our first example we will produce a finite dimensional example
with range not an arithmetic progression but such that $(\Cal L)$ uniquely
determines $\mu$.  (In the closing remarks in [N] it was claimed that no
such example exists.)  For convenience we will make our measure
integer valued.
 
\example{Example 1}  Let $\mu$ have nine atoms with measures 1, 2, 5, 6, 7,
8, 9, 10, 11.  With the notation above and $a_1=1,a_2=2,...,a_9=11
,$ the
set $\Cal R$ contains the following elements.\newline 
\centerline{(1, 0, 1,-1, 0, 0, 0, 0, 0)}
\centerline{(1, 0, 0, 1,-1, 0, 0, 0, 0)}
\centerline{(1, 0, 0, 0, 1,-1, 0, 0, 0)}
\centerline{(1, 0, 0, 0, 0, 1,-1, 0, 0)}
\centerline{(1, 0, 0, 0, 0, 0, 1,-1, 0)}
\centerline{(1, 0, 0, 0, 0, 0, 0, 1,-1)}
\centerline{(0, 1, 1, 0,-1, 0, 0, 0, 0)}
\centerline{(0, 0, 1, 1, 0, 0, 0, 0,-1)}
It is not hard to see that these are linearly independent and
therefore by Proposition 4.1 $\mu$ is uniquely determined by $(\Cal L)$.  Clearly
4 is omitted from the range of $\mu$ but 1, 2 and 3 are included, so the
range is not an arithmetic progression.
\endexample
Our next two examples show that properties of the range are not
good indicators of whether uniqueness holds.
 
\example{Example 2}  We begin with two finite dimensional examples.  Let
the five atoms of $\mu$ be of size 1, 2, 2, 2, 5 and let the four atoms of $
\mu^{\prime}$ be of 
size 1, 2,4, 5.  The range of $\mu$ is the same as the range $\mu^{
\prime}$, namely 
$\{1,$ $2,...,12\}$.  For $\mu$ note that the 
set $\Cal R$ contains the four linearly independent
elements (1, 1, 1, $0,-1)\text{, (0, 1,-1, 0, 0), (0, 0, 1,-1, 0), 
and (1,-1,-1,-1, 1).}$  Thus $
\mu$ 
is uniquely determined by $(\Cal L)$.  On the other hand for $\mu^{\prime}$ 
we have
$\Cal R=\{\pm (1,\text{ 0, 1,-1), }$ $\pm (1\text{,-1,-1, $1)\}$}$ and 
thus dim sp $
\Cal R=2<3$.  Thus $\mu^{\prime}$ is not 
uniquely determined by $(\Cal L)$.  (Take $\nu$ to have atoms of size 
1, 2, 6, and 7, 
for example.)\endexample
 
\example{Example 3}  We will modify the measures in Example 2 to get
two measures with range the same finite union of intervals, but so
that uniqueness holds for the first but not the second.  To do this
observe that if we omit the first atom from each of the measures in
Example 2, the range of the restricted measure in each case is
$\{\text{2, 4, 5, 6, 7, 9, $11$ $\}$}$.  Let $\gamma$ be any 
probability measure on a
measurable space $(\Omega_1,\Cal C)$ with range [0,1].  Replace the atom $
A_1$ by $\Omega_1$
and enlarge the $\sigma$-algebra to contain $\Cal C$ for each of $
\mu$ and $\mu^{\prime}$ and
define $\mu (A)=\mu^{\prime}(A)=\gamma (A)\text{ for }A\in \Cal C$.  (Extend $
\mu$ and $\mu^{\prime}$ to the
$\sigma$-algebra in the obvious way.)  Then the range of $\mu$ and $
\mu^{\prime}$ is now
$[0,12]\setminus ((1,2)\cup (3,4)\cup (8,9)\cup (10,11))$.  
Proposition 1.2 implies that if $
\nu$ is
another finite measure on $\Omega_1\cup\cup_{i>1}A_i$ satisfying $
(\Cal L)$ for either $\mu$ or $\mu^{\prime}$
then $\nu$ is a multiple of $\gamma$ on $\Omega_1$.  
Thus the uniqueness problem is the
same here as in Example 2, i.e., for $\mu$ condition $(\Cal L)$ implies uniqueness
but for $\mu^{\prime}$ it does not.
\endexample
These examples indicate that a characterization of the measures for
which $(\Cal L)$ is sufficient for uniqueness would require a similar result
for the finite dimensional case and thus seems to be difficult.
Finally we present an example of a finite signed measure with range
the interval [-1,1] which is uniquely determined by $(\Cal L)$ but for which
neither the positive nor negative parts are uniquely determined by
$(\Cal O)$. 
 
\example{Example 4}  For each natural number $n$ let $A_{2n-1}$ and $
A_{2n}$ be
atoms of a signed measure $\mu$ with $\mu (A_{2n-1})=\frac 2{3^n}=
-\mu (A_{2n})$.  It is not
hard to see that the ranges of $\mu^{+}$ and $\mu^{-}$ are both the usual Cantor
set $\Cal C$ in [0,1].  Moreover it is known and elementary to see that
$\Cal C-\Cal C=[-1,1].$ 
By Lemma 2.1, $|\mu|$ is uniquely determined by $(\Cal L)$. Indeed,
for each $n$, $1_{A_{2n-1}}-1_{\cup_{j>2n}A_j}$ and
$1_{A_{2n}}-1_{\cup_{j>2n}A_j} \in X$. Also
$1_{A_{2n-1}}-\frac121_{\cup_{j\geq 2n}A_j}=1_{A_{2n-1}}-1_{\cup_{j>2n}A_j}
-\frac12(1_{A_{2n}}-1_{\cup_{j>2n}A_j}) \in X$. Thus $X_\infty \supset
L^0_\infty(|\mu|)$. Lemma 1.1 implies that $\mu$ is also uniquely
determined by $(\Cal L)$.
  On the other hand
$\frac 2{3^n}>\frac 1{3^n}=\sum_{i>n}\frac 2{3^i}$ and thus for $\mu^{
+}$ and $\mu^{-}$ each atom is a bully. \endexample

Notice that there is no restriction imposed on $\dsize
\frac{d\nu}{d\mu}$ in this example and thus Theorem 3.5 does not
apply. It seems likely that the restriction is not necessary in
general.

\enddocument
\bye